# Geometry of 3-Dimensional Gradient Ricci Solitons with Positive Curvature

Sun-Chin Chu
National Chung Cheng University
Ming-Hsiung, Chia-Yi 621
Taiwan
scchu@math.ccu.edu.tw

## 1 Introduction

A solution $g(t)$ of the unnormalized Ricci flow on a differentiable manifold $\mathcal{M}$ is called a steady Ricci soliton if there exists a one-parameter family of diffeomorphisms $\varphi(t) : \mathcal{M} \to \mathcal{M}$ such that $g(t) = \varphi(t)^* g(0)$. Such solutions are fixed points of the unnormalized Ricci flow in the space of metrics modulo diffeomorphisms. Differentiating this equation with respect to time implies that $\frac{\partial}{\partial t} g = \mathcal{L}_V g$, where $V(t)$ is the one-parameter family of vector fields generated by $\varphi(t)$. By definition of the Ricci flow, the evolution equation of $g(t)$ is equivalent to $-2R_{ij} = \nabla_i V_j + \nabla_j V_i$. If $V$ is the negative gradient of a function $f$, then the solution is called a gradient Ricci soliton and the function $f$ is called a Ricci potential. On a gradient soliton, the Ricci tensor $R_{ij}$ equals $\nabla_i \nabla_j f$. Gradient Ricci solitons play a fundamental role in the study of singularities. In particular, they arise as singularity models (the limits of dilations about a singularity).

In dimension two, Hamilton's cigar is the only nonflat gradient Ricci soliton up to homothety. Under the Ricci flow, the metric evolves along the conformal vector field $V = \frac{\partial}{\partial r} = x\frac{\partial}{\partial x} + y\frac{\partial}{\partial y}$. More precisely, Hamilton's cigar $\Sigma^2$ is the complete surface with underlying manifold $\mathbb{R}^2$ and metric $g_\Sigma \doteq (dx^2 + dy^2)/(1 + x^2 + y^2)$. The cigar soliton is asymptotic to a cylinder of radius 1, while $R$ falls off like $e^{-s}$, where $s$ is the distance from some origin. R. Bryant and T. Ivey [I] show that there exists a complete, rotationally symmetric gradient soliton on $\mathbb{R}^n$, for $n \geq 3$, which is unique up to homothety. While this metric cannot be written down explicitly, one can compute that it has strictly positive sectional curvature that falls off like $1/s$.

Suppose that $(\mathcal{M}^n, g(t))$ is a gradient Ricci soliton with strictly positive curvature operator, whose scalar curvature assumes its maximum at some point $\mathcal{O}$, called the origin. Instead of considering the evolution equation of $g$, we focus attention on the level sets of a Ricci potential $f$. For convenience of exposition, let's recall some basic results on convex functions and their level sets. We say that a smooth function $\varphi$ is strictly convex if $\nabla\nabla\varphi$ is positive definite at each point, and $\varphi$ is convex if $\nabla\nabla\varphi$ is positive semidefinite, at each point. Bishop and O'Neill (Proposition 2.1 in [BO]) show that the critical points of a convex function are its absolute minimum points. So by strict convexity of $f$ and the fact that $R + |\nabla f|^2 = R(\mathcal{O})$ (see Theorem 20.1 in [H]), the critical point $\mathcal{O}$ is the only minimum point of $f$. Without loss of generality assume that the Ricci potential $f$ takes the value 0 at $\mathcal{O}$. Here and throughout, let $\mathcal{S}_\lambda$ be the level sets $f^{-1}(\lambda)$ and $\mathcal{R}_\lambda$ the sublevel sets $f^{-1}([0,\lambda])$. Remark 2.7 in [BO] says that if $\varphi$ is a convex function such that $\varphi^{-1}(c)$ is compact for some $c$, then for any two values $a, b > c$ of $\varphi$ the level sets $\varphi^{-1}(a)$ and $\varphi^{-1}(b)$ are compact and diffeomorphic under the flow of the vector field $\nabla\varphi/|\nabla\varphi|^2$. Hence all level sets of $f$ are compact, and the level sets $\mathcal{S}_{\lambda_1}$ and $\mathcal{S}_{\lambda_2}$ are diffeomorphic for all $\lambda_1, \lambda_2 > 0$.

In this article, we mainly study 3-dimensional complete gradient Ricci solitons with positive sectional curvature, whose scalar curvature attains its maximum at some point. Throughout this paper, we denote such solitons by $(\mathcal{M}^3, g)$. This paper is organized as follows. In section 2, we estimate the area growth of level sets and the volume growth of sublevel sets of $f$. These estimates are crucial in the study of the diameter of level sets and the geometry of such solitons at infinity. In section 3, we show that the scalar curvature of such solitons approaches zero at infinity. The vanishing of the scalar curvature at infinity plays a significant role in Hamilton's dimension reduction argument on odd-dimensional gradient solitons (see §22 [H]). In section 4, we first investigate the growth rate of the diameter of level sets, which implies that the tangent cone of such solitons is a ray $\mathbb{R}^+ \cup \{0\}$. Next we show that the scalar curvature falls off like $1/s$ and the diameter of geodesic spheres grows like $\sqrt{s}$ provided that $\mathbb{R} \times \Sigma^2$ cannot occur as a limit from dimension reduction. Note that the results in §4 [P] imply that $\mathbb{R} \times \Sigma^2$ cannot form as a limit of dilations of a compact solution. Therefore, if a gradient soliton $(\mathcal{M}^3, g)$ arises as a limit of dilations of a compact solution, then $\mathbb{R} \times \Sigma^2$ cannot occur as a limit from some dimension reduction since a limit of a limit is again a limit.

# 2  Volume growth of gradient Ricci solitons

Here and throughout, let $\mathcal{A}_\lambda$ denote the area of $\mathcal{S}_\lambda$, $\mathcal{V}_\lambda$ the volume of $\mathcal{R}_\lambda$, and $N$ the outward normal $\nabla f/|\nabla f|$ to the level sets. In the first part of this section we study the volume growth of odd-dimensional gradient Ricci solitons. In the second part of the section, we consider 3-dimensional gradient solitons and estimate the growth rates of $\mathcal{A}_\lambda$ and $\mathcal{V}_\lambda$ in terms of $\lambda$. Theorem 2.3 plays a key role in the study of the diameter of level sets and the geometry of such solitons at infinity in §4.

## 2.1  Volume growth is not large

Suppose that $(M^n, g)$ is a complete Riemannian $n$-manifold with nonnegative Ricci tensor. Let $B(p, r)$ be a metric ball of radius $r$ around $p$. The Bishop-Gromov theorem [GLP] says that the function

$$r \to \frac{\mathbf{vol}(B(p,r))}{\omega_n r^n}$$

is monotone decreasing for any $p \in M$. Define the asymptotic volume ratio $\alpha_M$ by

$$\alpha_M \doteqdot \lim_{r \to \infty} \frac{\mathbf{vol}(B(p,r))}{\omega_n r^n}.$$

It is known that
$$\alpha_M \omega_n r^n \leq \mathbf{vol}(B(p,r)) \leq \omega_n r^n,$$

$\alpha_M$ is independent of $p$, and is scale invariant.

**Definition**. $(M, g)$ is said to have large volume growth if $\alpha_M > 0$.

Let $(M_i, p_i)$ be a sequence of complete Riemannian $n$-manifolds with nonnegative Ricci tensor. By Gromov's precompactness theorem [GLP], $(M_i, p_i)$ subconverges to a complete length space $(Y, y)$. If $(M_i, p_i)$ satisfies the additional condition

$$\mathbf{vol}(B(p_i, 1)) \geq v > 0,$$

then it is known that (1) $Y$ has Hausdorff dimension $n$ and (2) for all $r > 0$, if $x_i \in M_i$ with $x_i \to x \in Y$, then

$$\mathbf{vol}(B(x_i, r)) \to \mathbf{vol}(B(x, r)) \text{ as } i \to \infty.$$

These results imply that the asymptotic volume ratios $\alpha_{M_i}$ and $\alpha_Y$ satisfy the inequality

$$\alpha_Y \geq \limsup_{i \to \infty} \alpha_{M_i}. \tag{2.1}$$

We have the following.

**Theorem 2.1** *Suppose that $(\mathcal{M}^{2n+1}, g)$ is an odd-dimensional gradient Ricci soliton with strictly positive curvature operator, whose scalar curvature assumes its maximum at some point. Assume that the non-flat factor from some dimension reduction does not have large volume growth. Then the soliton does not have large volume growth. In particular, if the non-flat factor is compact, then $\mathcal{M}$ cannot have large volume growth.*

**Proof.** Since the asymptotic scalar curvature ratio $\limsup_{s \to \infty} Rs^2 = \infty$ and $\dim \mathcal{M} = 2n+1$ is odd, we can perform dimension reduction and have a limit which splits as a product with a line $\mathbb{R}$ of an ancient solution of bounded positive curvature. (See Theorem 22.3 in [H].) By hypothesis, the asymptotic volume ratio of the limit manifold is zero. As a result of (2.1), we have $\alpha_\mathcal{M} = 0$. ∎

**Corollary 2.2** *Suppose that $(\mathcal{M}^3, g)$ is a gradient Ricci soliton with positive sectional curvature, whose scalar curvature assumes its maximum at some point. Then $(\mathcal{M}^3, g)$ cannot have large volume growth.*

**Proof.** We can classify the non-flat surface as a round sphere $\mathbb{S}^2$ or the cigar soliton $\Sigma^2$. It is clear that the volume growth of either $\mathbb{R} \times \mathbb{S}^2$ or $\mathbb{R} \times \Sigma^2$ is not large, so the volume growth of $(\mathcal{M}^3, g)$ cannot be large. ∎

## 2.2 Growth rates of the area of level sets and volume of the sublevel sets

Let $e_1, \ldots, e_{n-1}$ be an orthonormal frame on $\mathcal{S}_\lambda$. Denote by $K_{ij} = K_\mathcal{M}(e_i, e_j)$ the sectional curvature of the plane spanned by the vectors $e_i, e_j$, and by $Ric(v)$ the Ricci tensor in the direction $v$. It is known that the second fundamental form $II_\lambda$ on $\mathcal{S}_\lambda$ is given by

$$II_\lambda(e_i, e_j) = \langle \nabla_{e_i} N, e_j \rangle = \frac{Ric(e_i, e_j)}{|\nabla f|}.$$

The following theorem plays an important role in the study of 3-dimensional gradient Ricci solitons with positive sectional curvature.

**Theorem 2.3** *Suppose that $(\mathcal{M}^3, g)$ is a complete gradient Ricci soliton with positive sectional curvature, whose scalar curvature assumes its maximum at some point $\mathcal{O}$. Then we have the following:*
*(1) The area of level sets $\mathcal{A}_\lambda$ grows linearly in $\lambda$.*
*(2) The volume of sublevel sets $\mathcal{V}_\lambda$ grows quadratically in $\lambda$.*

**Proof.** Without loss of generality, we may assume that the maximum scalar curvature is 1. Let $\Gamma(\lambda)$ denote the maximal integral curve of the vector field $\nabla f / |\nabla f|^2$. From the fact that

$$\frac{d}{d\lambda} f(\Gamma(\lambda)) = \nabla f(\Gamma(\lambda)) \cdot \dot{\Gamma}(\lambda) = \nabla f \cdot \frac{\nabla f}{|\nabla f|^2} = 1,$$

we see that $f$ takes the value $\lambda$ at $\Gamma(\lambda)$ and the flow generated by the vector field $\nabla f / |\nabla f|^2$ evolves level sets of $f$ into themselves. Evolving level sets $\mathcal{S}_\lambda$ along the vector field $\nabla f / |\nabla f|^2$, the area $\mathcal{A}_\lambda$ changes at a rate

$$\frac{d\mathcal{A}_\lambda}{d\lambda} = \int_{\mathcal{S}_\lambda} HN \cdot \frac{\nabla f}{|\nabla f|^2} \, da,$$

where $H$ stands for the mean curvature of level sets $\mathcal{S}_\lambda$. By straightforward computation,

$$\frac{d\mathcal{A}}{d\lambda} = \int_{\mathcal{S}_\lambda} HN \cdot \frac{\nabla f}{|\nabla f|^2} da = \int_{\mathcal{S}_\lambda} \frac{Ric(e_1) + Ric(e_2)}{|\nabla f|^2} > 2 \int_{\mathcal{S}_\lambda} K_{\mathcal{M}}(e_1, e_2) da$$

since $|\nabla f| < 1$, $Ric(e_1) + Ric(e_2) = 2K_{\mathcal{M}}(e_1, e_2) + Ric(N)$, and $Ric(N) > 0$.

**Claim.** *There exist constants $c, \lambda(c) > 0$ such that*

$$c < \int_{\mathcal{S}_\lambda} K_{\mathcal{M}}(e_1, e_2) da < 4\pi \tag{2.2}$$

*for $\lambda \geq \lambda(c)$.*
**Proof.** We will show that the integral in (2.2) is monotone increasing.

**Lemma 2.4** *Suppose that $(\mathcal{M}^n, g)$, for $n \geq 3$, is a gradient Ricci soliton with positive sectional curvature, whose scalar curvature assumes its maximum at some point $\mathcal{O}$. Then we have*

$$\frac{d}{d\lambda} \int_{\mathcal{S}_\lambda} \det(II_\lambda) da < 0. \tag{2.3}$$

**Proof.** Define a $\wedge^n T\mathcal{M}$-valued $(n-1)$-form $\alpha$ by

$$\alpha(X_1, \ldots, X_{n-1}) = N \wedge \nabla_{X_1} N \wedge \cdots \wedge \nabla_{X_{n-1}} N$$

for vectors $X_1, \ldots, X_{n-1}$. Evaluate $\alpha$ at an orthonormal frame $e_1, \ldots, e_{n-1}$ as follows.

$$\alpha(e_1, \ldots, e_{n-1}) = N \wedge \nabla_{e_1} N \wedge \cdots \wedge \nabla_{e_{n-1}} N = \det(II_\lambda) N \wedge e_1 \wedge \cdots \wedge e_{n-1}.$$

By identifying the bundle $\wedge^n T\mathcal{M}$ with $\mathcal{M} \times \mathbb{R}$, we have

$$\int_{\mathcal{S}_\lambda} \det(II_\lambda) = \int_{\mathcal{S}_\lambda} \alpha.$$

Denote by $\mathcal{R}_{\lambda_1}^{\lambda_2}$ the region $f^{-1}([\lambda_1, \lambda_2])$ for $\lambda_1 < \lambda_2$. By the divergence theorem,

$$\int_{\mathcal{S}_{\lambda_2}} \det(II_{\lambda_2}) - \int_{\mathcal{S}_{\lambda_1}} \det(II_{\lambda_1}) = \int_{\mathcal{R}_{\lambda_1}^{\lambda_2}} d^\nabla \alpha. \tag{2.4}$$

Compute the exterior differential of the $\wedge^n T\mathcal{M}$-valued $(n-1)$-form $\alpha$ as follows. Evaluating $d^\nabla \alpha$ at an orthonormal frame $e_0 = N, e_1, \ldots, e_{n-1}$,

$$\begin{aligned}
&(d^\nabla \alpha)(e_0, e_1, \ldots, e_{n-1}) \\
&= \sum_i (-1)^i \nabla_{e_i}(\alpha(N, \ldots, \hat{e}_i, \ldots, e_{n-1})) \\
&\quad + \sum_{i \neq j} (-1)^{i+j} \alpha([e_i, e_j], N, \ldots, \hat{e}_i, \ldots, \hat{e}_j, \ldots, e_{n-1}).
\end{aligned}$$

Since $\nabla_N N = 0$, by straightforward computation we have

$$\begin{aligned}
&(d^\nabla \alpha)(N, e_1, \ldots, e_{n-1}) \\
&= \nabla_N(\alpha(e_1, \ldots, e_{n-1})) + \sum_{i \neq j} (-1)^j \alpha([N, e_j], e_1, \ldots, \hat{e}_j, \ldots, e_{n-1}) \\
&= N \wedge (\nabla_N \nabla_{e_1} N) \wedge \cdots \wedge \nabla_{e_{n-1}} N + \cdots + N \wedge \cdots \wedge \nabla_{e_{n-2}} N \wedge (\nabla_N \nabla_{e_{n-1}} N) \\
&\quad + N \wedge \nabla_{[N, e_1]} N \wedge \cdots \wedge \nabla_{e_{n-1}} N + \cdots + N \wedge \cdots \wedge \nabla_{e_{n-2}} N \wedge \nabla_{[N, e_{n-1}]} N \\
&= N \wedge (\nabla_N \nabla_{e_1} N - \nabla_{e_1} \nabla_N N - \nabla_{[N, e_1]} N) \wedge \nabla_{e_2} N \wedge \cdots \wedge \nabla_{e_{n-1}} N + \cdots \\
&\quad + N \wedge \cdots \wedge \nabla_{e_{n-2}} N \wedge (\nabla_N \nabla_{e_{n-1}} N - \nabla_{e_{n-1}} \nabla_N N - \nabla_{[N, e_{n-1}]} N) \\
&= N \wedge (R(N, e_1) N) \wedge \cdots \wedge \nabla_{e_{n-1}} N + \cdots + N \wedge \cdots \wedge \nabla_{e_{n-2}} N \wedge (R(N, e_{n-1}) N) \\
&= -[N \wedge (R(e_1, N) N) \wedge \cdots \wedge \nabla_{e_{n-1}} N + \cdots + N \wedge \cdots \wedge \nabla_{e_{n-2}} N \wedge (R(e_{n-1}, N) N)].
\end{aligned}$$

By choosing the orthonormal principal directions as the $(n-1)$ vectors $e_1, \ldots, e_{n-1}$, we have

$$(d^\nabla \alpha)(e_0, e_1, \ldots, e_{n-1}) = -(n-2)! \Sigma_{j=1}^{n-1} \left[ K_{0j} \left( \Pi_{i>0,\ i \neq j}\ \mu_i \right) \right],$$

where $K_{0j} = \langle R(e_j, N)N, e_j \rangle$ and $II_\lambda(e_i) = \mu_i e_i$.

Since $K_{0j} > 0$ and $\mu_i > 0$ for all $1 \leq i, j \leq n-1$, the right-hand side of (2.4) is negative, thus

$$\int_{\mathcal{S}_{\lambda_2}} \det(II_{\lambda_2}) < \int_{\mathcal{S}_{\lambda_1}} \det(II_{\lambda_1}).$$

The lemma follows. ■

Let $\mathcal{K}_\lambda$ denote the sectional curvature of the induced metric on $\mathcal{S}_\lambda$. The Gauss-Bonnet theorem implies that

$$\int_{\mathcal{S}_\lambda} K_{\mathcal{M}}(e_1, e_2) = \int_{\mathcal{S}_\lambda} \mathcal{K}_\lambda - \int_{\mathcal{S}_\lambda} \det(II_\lambda) = 4\pi - \int_{\mathcal{S}_\lambda} \det(II_\lambda)$$

is strictly increasing by the monotonicity of (2.3). Consequently, the claim follows. ■

By the coarea formula, we have the following equations:

$$\frac{d}{d\lambda} \mathcal{V}_\lambda = \frac{d}{d\lambda} \int_{\mathcal{R}_\lambda} dv = \int_{\mathcal{S}_\lambda} \frac{1}{|\nabla f|},$$

$$\frac{d}{d\lambda} \int_{\mathcal{S}_\lambda} \frac{1}{|\nabla f|} = \int_{\mathcal{S}_\lambda} \frac{1}{|\nabla f|} \cdot \frac{1}{|\nabla f|} \cdot \frac{\sum_{i=1}^{n-1} Ric(e_i)}{|\nabla f|} - \frac{\nabla_i \nabla_j f \cdot \nabla_i f}{|\nabla f|^3} \cdot \frac{\nabla_j f}{|\nabla f|^2}$$

$$= \int_{\mathcal{S}_\lambda} \frac{\sum_{i=1}^{n-1} \sum_{j=0,\ j \neq i}^{n-1} K_{ij} - \sum_{1 \leq i \leq n-1} K_{i0}}{|\nabla f|^3}$$

$$= 2 \int_{\mathcal{S}_\lambda} \frac{\sum_{1 \leq i < j \leq n-1} K_{ij}}{|\nabla f|^3} > 0.$$

Therefore, $\mathcal{V}_\lambda$ is strictly convex as a function of $\lambda$. In particular, for a three-manifold we have

$$\frac{d}{d\lambda} \int_{\mathcal{S}_\lambda} \frac{1}{|\nabla f|} = 2 \int_{\mathcal{S}_\lambda} \frac{K_{\mathcal{M}}}{|\nabla f|^3}.$$

Since $R + |\nabla f|^2 = 1$ and $\nabla_i \nabla_j f = R_{ij} > 0$, we have

$$\frac{d}{d\lambda} R(\Gamma(\lambda)) = \nabla R \cdot \dot{\Gamma}(\lambda) = -2 R_{ij} \nabla_i f \cdot \dot{\Gamma}_j(\lambda)$$

$$= -\frac{2 \nabla_i \nabla_j f \cdot \nabla_i f \cdot \nabla_j f}{|\nabla f|^2} < 0, \tag{2.5}$$

which means that the scalar curvature $R$ decreases along the integral curves. Note that this also implies that $|\nabla f|$ increases along the integral curves. Thus, there exist constants $\delta, \lambda(\delta) > 0$ such that $1 - \delta \leq |\nabla f|^3 < 1$ and

$$2\int_{\mathcal{S}_\lambda} K_\mathcal{M} dA < \frac{d}{d\lambda}\int_{\mathcal{S}_\lambda} \frac{1}{|\nabla f|} dA \leq \frac{2}{1-\delta}\int_{\mathcal{S}_\lambda} K_\mathcal{M} dA \qquad (2.6)$$

for $\lambda \geq \lambda(\delta)$. (Note that we shall prove that the scalar curvature goes to zero at infinity under the same hypothesis as in this theorem. This means that we can choose $\delta > 0$ arbitrarily.)

The area estimate follows from the estimates (2.2), (2.6) and the fact that

$$\mathcal{A}_\lambda < \int_{\mathcal{S}_\lambda} \frac{1}{|\nabla f|} \leq \frac{\mathcal{A}_\lambda}{1-\delta}.$$

Similarly, by the same estimates (2.2), (2.6) and the identity

$$\mathbf{vol}(\mathcal{R}_{\lambda_0}^\lambda) = \int_{\lambda_0}^{\lambda} \left(\int_{\mathcal{S}_r} \frac{1}{|\nabla f|}\right),$$

there exist constants $c', c'', \lambda_0 > 0$ such that

$$c''(\lambda^2 - \lambda_0^2) \leq \mathbf{vol}(\mathcal{R}_{\lambda_0}^\lambda) \leq c'(\lambda^2 - \lambda_0^2) \qquad (2.7)$$

for $\lambda \geq \lambda_0$.

The volume estimate follows from the last estimate (2.7) and the relation $\mathbf{vol}(\mathcal{R}_\lambda) = \mathbf{vol}(\mathcal{R}_{\lambda_0}) + \mathbf{vol}(\mathcal{R}_{\lambda_0}^\lambda)$. ∎

By inequalities (2.2) and $R > 2K_\mathcal{M}(e_1, e_2)$ there exist constants $c, \lambda(c) > 0$ such that

$$\left(\max_{\mathcal{S}_\lambda} R\right) \mathcal{A}_\lambda \geq \int_{\mathcal{S}_\lambda} R > 2\int_{\mathcal{S}_\lambda} K_\mathcal{M} > 2c$$

for $\lambda \geq \lambda(c)$. We have the following.

**Corollary 2.5** *Under the same hypothesis we have*

$$\limsup_{d(p,\mathcal{O})\to\infty} R\lambda > 0.$$

# 3 Scalar curvature of 3-dimensional solitons with positive curvature

In this section, we first give an estimate for the angle $\theta$ between the outward normal $N$ and the radial direction $\partial_r$ at any point $q \in \mathcal{S}_\lambda$. Note that in the proof of Theorem 22.3 [H], Hamilton wrote that

> "... $|Df|^2$ approaches the maximum curvature $M$ as $s \to \infty$ where $s$ is the distance from some origin."

This is equivalent to $R \to 0$ at infinity. To the best of the author's knowledge, this has not been justified in any literature yet. Next we show that the scalar curvature goes to zero at infinity on 3-dimensional gradient Ricci solitons with positive sectional curvature and scalar curvature assuming its maximum at some point.

## 3.1 Ricci potential and distance function

It has been seen in section 2 that $|\nabla f|$ increases along the integral curves $\Gamma$. Let $\zeta_\Gamma$ denote the limit $\zeta_\Gamma \doteq \lim_{\lambda \to \infty} |\nabla f(\Gamma(\lambda))|$. To prove Theorem 3.5, we need the following.

**Theorem 3.1** *Suppose that $(\mathcal{M}^n, g)$ is a gradient Ricci soliton with positive sectional curvature, assuming its maximum scalar curvature 1 at some point $\mathcal{O}$. Assume that $\gamma(t)$ is an arbitrary minimizing geodesic from $\mathcal{O}$ to $q = \Gamma(\lambda) \in \mathcal{S}_\lambda$. Then given $\varepsilon > 0$ there exists a constant $\lambda(\varepsilon) > 0$ such that the angle $\theta(\Gamma, \gamma)$ between the outward normal $N$ and the radial direction $\partial_r$ at $\gamma(s) = q$ satisfies*

$$|\nabla f(q)| \cos \theta(\Gamma, \gamma) \geq \frac{f(q) - f(\mathcal{O})}{s}.$$

*Therefore, we have*

$$\cos \theta(\Gamma, \gamma) \geq (\zeta_\Gamma - \varepsilon)/\zeta_\Gamma \tag{3.1}$$

*for all $\lambda \geq \lambda(\varepsilon)$.*

**Proof.** Since the Ricci potential is geodesically convex, it suffices to find a lower bound for the Ricci potential in terms of the distance function.

**Lemma 3.2** *Under the previous hypothesis, given $\epsilon \in (0, \zeta_\Gamma)$ there is a constant $\lambda(\epsilon)$ such that*

$$(\zeta_\Gamma - \epsilon)s \leq f(\Gamma(\lambda)) < s \tag{3.2}$$

*for all $\lambda \geq \lambda(\epsilon)$.*

**Proof.** The right-hand side of inequality (3.2) follows easily from the mean value theorem and the fact that $|\nabla f| < 1$ on the soliton. We only need to prove the left-hand side of inequality (3.2).

For any given $\epsilon \in (0, \zeta_\Gamma)$, there is a number $\lambda_0$ such that $|\nabla f(\Gamma(\lambda_0))| = \zeta_\Gamma - \epsilon/2$ and $|\nabla f(\Gamma(\lambda))| > \zeta_\Gamma - \epsilon/2$ for all $\lambda > \lambda_0$. Let $\bar{\lambda}_0$ be any real number greater than $\lambda_0$ and define

$$\eta \doteq \inf_{\lambda \geq \bar{\lambda}_0} |\nabla f(\Gamma(\lambda))| = |\nabla f(\Gamma(\bar{\lambda}_0))| > \zeta_\Gamma - \epsilon/2.$$

The length of $\Gamma|_{[\bar{\lambda}_0, \lambda]}$ is given by

$$\int_{\bar{\lambda}_0}^{\lambda} |\Gamma'(\tau)| d\tau = \int_{\bar{\lambda}_0}^{\lambda} \frac{1}{|\nabla f(\Gamma(\tau))|} d\tau < \eta^{-1}(\lambda - \bar{\lambda}_0).$$

The distance from $\Gamma(\lambda)$ to the level set $\mathcal{S}_{\lambda_0}$ is at most the length of $\Gamma|_{[\bar{\lambda}_0, \lambda]}$. This implies that

$$d(\Gamma(\lambda), \mathcal{S}_{\bar{\lambda}_0}) < \eta^{-1}(\lambda - \bar{\lambda}_0).$$

Thus, we have

$$f(\Gamma(\lambda)) > \bar{\lambda}_0 + \eta d(\Gamma(\lambda), \mathcal{S}_{\bar{\lambda}_0}) > \bar{\lambda}_0 + (\zeta - \epsilon/2) d(\Gamma(\lambda), \mathcal{S}_{\bar{\lambda}_0}).$$

The left-hand side of inequality (3.2) follows from the last estimate and the fact that $f$ is an exhaustion function. ∎

Since $f$ is geodesically convex, we have

$$|\nabla f(q)||\gamma'(s)| \cos \theta(\Gamma, \gamma) = \frac{d}{dt}\Big|_{t=s} f(\gamma(t)) \geq \frac{f(q) - f(\mathcal{O})}{s}.$$

The theorem follows. ∎

As a consequence of Lemma 3.2, we have the following.

**Corollary 3.3** *Suppose that $(\mathcal{M}^n, g)$ is a complete gradient Ricci soliton with positive sectional curvature, and scalar curvature assuming its maximum 1 at some point $\mathcal{O}$ and going to zero at infinity. For any given $\delta > 0$, there exists an $\bar{s} = \bar{s}(\delta) > 0$ such that*

$$(1 - \delta)s \leq \lambda < s, \tag{3.3}$$

*where $\lambda = f(p)$ and $s = d(\mathcal{O}, p) \geq \bar{s}$ for $p \in \mathcal{M}$.*

*Moreover, the angle between the outward normal and the radial direction approaches zero uniformly at infinity.*

The following observation plays an important role in the study of the behavior of scalar curvature at infinity in §3.2.

**Theorem 3.4** *The affine space of Ricci potentials with linear growth on the product manifold $(\mathbb{R} \times \Sigma^2, ds^2 + d\rho^2 + \tanh^2 \rho \, d\theta^2)$ is given by*

$$\{c_1 + c_2 s + 2 \ln \cosh \rho \mid c_1, c_2 \in \mathbb{R}\},$$

*where $2 \ln \cosh \rho$ is a Ricci potential of Hamilton's cigar $(\Sigma^2, g_\Sigma)$.*

**Proof.** Let $f_1, f_2$ be two Ricci potentials with linear growth. The difference function $f_1 - f_2$ satisfies

$$\nabla^2 (f_1 - f_2) = \nabla^2 f_1 - \nabla^2 f_2 = Ric - Ric = 0;$$

thus $\nabla(f_1 - f_2)$ is parallel. If $f_1 - f_2$ is not of the form $c_1 + c_2 s$, then $\mathbb{R} \times \Sigma$ must have a line other than $\mathbb{R} \times \{p\}$ for all $p \in \Sigma$, which is impossible. The theorem follows. ∎

**Remark**. (1) In polar coordinates, the cigar metric is given by $g_\Sigma = (dr^2 + r^2 d\theta^2)/(1+r^2)$. Making the change of variables $\rho = \sinh^{-1} r = \ln(\sqrt{1+r^2} + r)$, we have $g_\Sigma = d\rho^2 + \tanh^2 \rho \, d\theta^2$. (2) Suppose that $g \doteq \frac{1}{a^2} g_\Sigma$, i.e.,

$$g \doteq \frac{1}{a^2}(d\rho^2 + \tanh^2 \rho \, d\theta^2) = d\sigma^2 + \frac{1}{a^2} \tanh^2(a\sigma) d\theta^2,$$

where $\sigma = \rho/a$ for some constant $a > 0$. Then the function

$$\pi(\sigma) = 2 \ln \cosh(a\sigma) = 2 \ln \cosh \rho$$

is a Ricci potential of $g$.

### 3.2 Scalar curvature vanishes at infinity

Since the scalar curvature $R$ decreases along the integral curves $\Gamma$, denote by $R_\Gamma$ the limit of $R$ along $\Gamma$. It is clear that we have either (i) $R_\Gamma = 0$ for all $\Gamma$ or (ii) $\hat{R} \doteq \sup_\Gamma R_\Gamma > 0$. In case (ii), we have either (ii-a) $\inf_\Gamma R_\Gamma > 0$ or (ii-b) $\inf_\Gamma R_\Gamma = 0$. The following theorem says that neither (ii-a) nor (ii-b) can occur.

**Theorem 3.5** *Suppose that $(\mathcal{M}^3, g)$ is a gradient Ricci soliton with positive sectional curvature and scalar curvature assuming its maximum at an origin $\mathcal{O}$. Then the scalar curvature goes to zero at infinity.*

**Proof.** We may assume that the maximum scalar curvature is 1. Suppose that we are in case (ii). From compactness of the collection of all maximal integral curves and the fact that the scalar curvature $R$ is monotone decreasing along the integral curves, there exists a maximal integral curve $\Gamma$ such that $R_\Gamma = \hat{R}$. To perform dimension reduction we need control over the injectivity radius. However, because in this case the scalar curvature $R$ decreases to $\hat{R} > 0$ along $\Gamma$, we do not need to scale the metric much when performing dimension reduction. Therefore, we get the required control of injectivity radius.

**Case (ii-a).** From the facts that the maximal integral curves of $\nabla f/|\nabla f|^2$ intersect the level sets $\mathcal{S}_\lambda$, $\lambda > 0$, the scalar curvature decreases along these curves, and $\inf_\Gamma R_\Gamma > 0$ is positive, the limit space must be $\mathbb{R} \times \mathbb{S}^2$ since the scalar curvature of $\mathbb{R} \times \Sigma^2$ goes to zero at infinity. According to the dimension reduction argument, we can choose a sequence of points $\{p_k\}$ on the curve $\Gamma$ such that the sequence of the pointed Riemannian manifolds $(\mathcal{M}^3, p_k, R(p_k)g)$ converges to the product manifold $\mathbb{R} \times \mathbb{S}^2$. We choose $k$ large enough so that there is a long topological and geometric neck $\mathcal{N}_k = (-L, L) \times \mathbb{S}^2$ with $p_k$ lying in the central slice $\mathcal{S}_0 = \{0\} \times \mathbb{S}^2$. Evolving the metric by $-2R_{ij}$ is equivalent deforming the slice $\mathcal{S}_\tau$ along the negative gradient vector field $-\nabla f$. Since $|\nabla f| < 1$, the slice $\mathcal{S}_\tau$ stays inside the long neck $\mathcal{N}_k$ for $\tau < \frac{L}{2}$ and the area of $\mathcal{S}_\tau$ is at least (roughly) $4\pi$. On the other hand, we have

$$\frac{dA_\mathcal{S}}{d\tau} = -\int_{\mathcal{S}_\tau} (2Rm(T) + Rc(N))\, da \leq -\frac{1}{2}\int_{\mathcal{S}_\tau} R\, da \leq -\frac{1}{2}\int_{\mathcal{S}_\tau} R\, da \approx -2\pi.$$

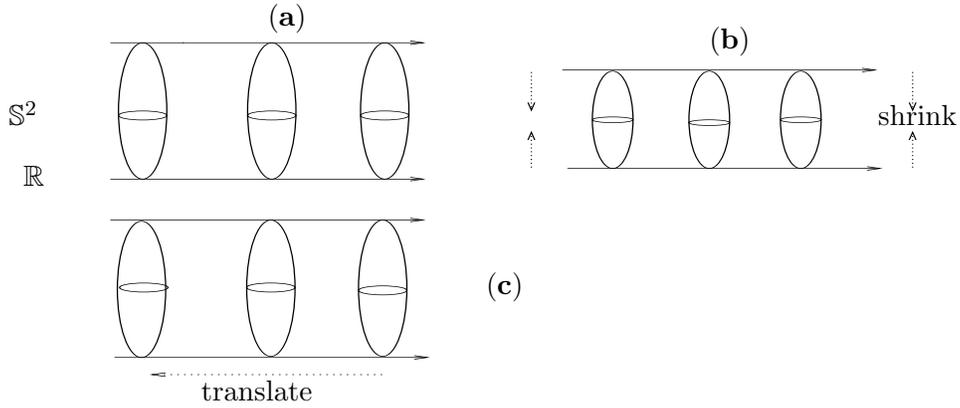

Figure 1: slices — $\mathbb{R} \times \mathbb{S}^2$

This implies that $area(\mathcal{S}_\tau)$ vanishes at a finite time $T \approx area(\mathcal{S})/(2\pi)$. This contradicts the previous observation, therefore, (ii-a) cannot occur.

**Case (ii-b).** Conversely, from the fact that $\inf_\Gamma R_\Gamma$ is zero the limit arising from the dimension reduction argument must be $\mathbb{R} \times \Sigma^2$ and there is a sequence of points $\{p_k\}$ on $\Gamma$ such that the sequence of pointed Riemannian manifolds $(\mathcal{M}, p_k, R(p_k)g)$ converges to the product manifold $\mathbb{R} \times \Sigma^2$ with metric

$$ds^2 + \frac{4}{\hat{R}}(d\rho^2 + \tanh^2\rho \, d\theta^2) = ds^2 + d\sigma^2 + \frac{4}{\hat{R}} \tanh^2(\sqrt{\hat{R}}\sigma/2) d\theta^2,$$

where $\sigma = 2\rho/\sqrt{\hat{R}}$ and $\rho = \sinh^{-1} r = \ln(\sqrt{1+r^2} + r)$.

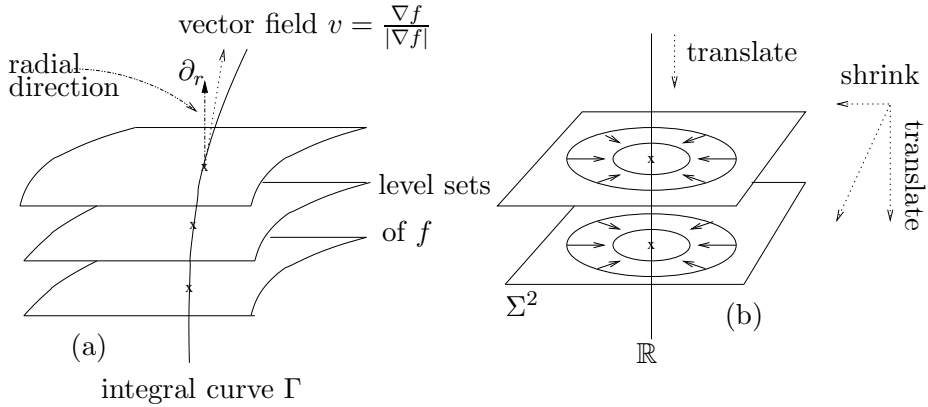

Figure 2: level sets — $\mathbb{R} \times \Sigma^2$

Since $(f(\Gamma(\lambda)))' = \nabla f \cdot \Gamma'(\lambda) = \nabla f \cdot \nabla f/|\nabla f|^2 = 1$ along the integral curve $\Gamma$, we have $f(\Gamma(\lambda)) = \lambda$ for all $\lambda \in [0, \infty)$ and $p_k = \Gamma(\lambda_k)$ for some $\lambda_k$. By local derivative estimates, $R_{ij} = \nabla_i \nabla_j (f - \lambda_k)$, and Theorem 3.4, we have that

$$f - f(p_k) \to \tilde{f}(s, \sigma, \theta) = \sqrt{1 - \hat{R}} s + \bar{f}(\sigma, \theta) \text{ as } k \to \infty$$

and

$$\tilde{R}(s, \sigma, \theta) + |\tilde{\nabla}\tilde{f}|^2 = 1,$$

where $\bar{f}(\sigma, \theta) = \bar{f}(\sigma) = 2\ln\cosh(\sqrt{\hat{R}}\sigma/2)$ is the Ricci potential of the cigar soliton $\Sigma^2$ with $\bar{R}(\mathcal{O}) = \hat{R}$. Note that $\bar{f}(\sigma)$ approaches $\sqrt{\hat{R}}\sigma$ as $\sigma \to \infty$.

As shown in Figure 3, the geodesic intersects the $s$-axis at a large angle, which contradicts the angle estimate (3.1) in Theorem 3.1. Therefore, case (ii-b) can not occur. ∎

As a consequence of Theorems 2.3, 3.5 and Proposition 3.3, we have that (1) the area of $\partial(B(\mathcal{O}, s))$ grows linearly and the volume of $B(\mathcal{O}, s)$ grows quadratically in $s$, and (2) $\limsup\limits_{s\to\infty} Rs > 0$.

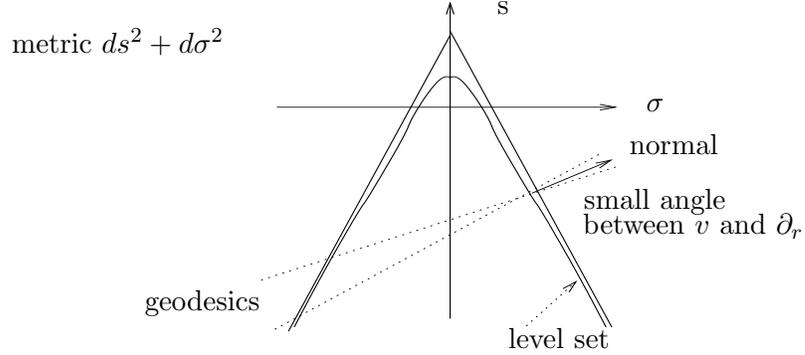

Figure 3: level set $\tilde{f}^{-1}(c) \cap \{\theta = 0\}$

## 4 Geometry of gradient Ricci soliton at infinity

We first recall the basic facts about the ideal boundary and the tangent cone of nonnegatively curved manifolds. As a consequence of Theorem 4.2, the tangent cone of $(\mathcal{M}^3, g)$ is a ray. Next we prove Theorem 4.3 by choosing the reference points $p_k$ of a sequence of pointed Riemannian manifolds suitably and scaling the metric by factors other than $R(p_k)$. Theorems 4.2 and 4.3 are important in the study of the asymptotic behavior of 3-dimensional gradient Ricci solitons arising as a limit of dilations of a compact solution.

### 4.1 Ideal boundary and tangent cone

Let $(\mathbb{M}, g)$ be a complete, noncompact Riemannian manifold with nonnegative sectional curvature. The ideal boundary of $(\mathbb{M}, g)$ consists of equivalence classes of rays $\mathcal{R}_o$ emanating from a base point $o$. Two rays $\sigma$ and $\gamma$ are said to be equivalent if $\lim_{t\to\infty} d(\sigma(t), \gamma(t))/t = 0$. The distance between two equivalence classes $[\sigma], [\gamma]$ is defined to be

$$d_\infty([\sigma], [\gamma]) = \lim_{t \to \infty} \frac{d_t(\sigma(t), \gamma(t))}{t},$$

where $d_t$ is the inner metric on geodesic sphere $S(o, t)$ of radius $t$ around the base point $o$. Denote by $(\mathbb{M}_\infty, d_\infty)$ the ideal boundary of $(\mathbb{M}, g)$. The following results related to the ideal boundary are useful for the rest of this paper. For manifolds of nonnegative sectional curvature, it is known that the Gromov-Hausdorff limit of $(\mathbb{M}, o, g/\lambda)$ as $\lambda \to \infty$ exists and is a cone over a compact metric space, which is isometric to the ideal boundary $\mathbb{M}_\infty$. This metric cone is known as the tangent cone of the nonnegatively curved manifold $(\mathbb{M}, g)$. Kasue (Proposition 2.3 in [K]) proves that

$$\lim_{t \to \infty} (S(o, t), d_t/t) = (\mathbb{M}_\infty, d_\infty). \tag{4.1}$$

By the convergence (4.1) and the estimate (3.3) the sequence of rescaled level sets $(\mathcal{S}_\lambda, d_\lambda/\lambda)$ converges to the ideal boundary $(\mathbb{M}_\infty, d_\infty)$ as $\lambda \to \infty$, where $d_\lambda$ is the inner metric on $\mathcal{S}_\lambda$, provided that $(\mathbb{M}, g)$ is a gradient soliton of positive curvature operator, whose scalar curvature goes to zero at infinity.

Suppose that the sectional curvature of $(\mathbb{M}, g)$ is bounded from above by a positive constant. Proposition 5.5 in [K] says that for each end $\mathcal{E}^\alpha$ of $\mathbb{M}$, one has

$$\liminf_{t \to \infty} \frac{\log \mathbf{vol}(B(p, t) \cap \mathcal{E}^\alpha)}{\log t} \geq 1 + \dim_\mathcal{H} \mathbb{M}_\infty^\alpha,$$

where $p$ is a fixed point of $\mathbb{M}$ and $\dim_\mathcal{H} \mathbb{M}_\infty^\alpha$ stands for the Hausdorff dimension of the connected component $\mathbb{M}_\infty^\alpha$ of $\mathbb{M}_\infty$ corresponding to $\mathcal{E}^\alpha$. From the volume estimate of geodesic balls, we have the following.

**Proposition 4.1** *Suppose that $(\mathcal{M}^3, g)$ is a gradient Ricci soliton with positive sectional curvature and scalar curvature assuming its maximum at some point. Then the Hausdorff dimension of the ideal boundary $\mathcal{M}_\infty$ is at most 1.*

### 4.2 Diameter of level sets

Greene and Shiohama (Proposition 1.4 in [GS]) show that the diameter $\mathcal{D}_\lambda$ of level sets is a nondecreasing function of $\lambda$. The rest of this paper is devoted to the study of the ratio $\mathcal{D}_\lambda/\lambda$ and the boundedness of the quantity $R\mathcal{D}_\lambda^2$. Corollary 4.5 gives the lower and upper bounds on the scalar curvature and the diameter $\mathcal{D}_\lambda$ provided that $\mathbb{R} \times \Sigma^2$ cannot occur as a limit of dimension reduction.

**Theorem 4.2** *Suppose that $(\mathcal{M}^3, g)$ is a gradient Ricci soliton with positive sectional curvature and scalar curvature assuming its maximum at an origin $\mathcal{O}$. Then we have*

$$\lim_{\lambda \to \infty} \frac{\mathcal{D}_\lambda}{\lambda} = 0. \tag{4.2}$$

*Therefore, the tangent cone is a ray.*

**Proof.** Grove and Petersen (Theorem A and Remark 2.4 in [GP]) show that a space arising as a Gromov-Hausdorff limit of Riemannian manifolds of fixed dimension with curvature bounded uniformly below has constant dimension.

The Gromov-Hausdorff limit of a sequence of such $n$-dimensional manifolds is either again of dimension $n$, perhaps with some metric singularities, or alternatively collapse occurs everywhere, and the limit space has everywhere a fixed dimension $k < n$. So it follows from Proposition 4.1 that we have either (i) $\dim_{\mathcal{H}} \mathcal{M}_\infty = 1$ or (ii) $\dim_{\mathcal{H}} \mathcal{M}_\infty = 0$. We claim that case (i) cannot occur as follows.

Since the area and the inner diameter of level sets $\mathcal{S}_\lambda$ grow linearly in $\lambda$, by Croke's estimate (Corollary 15 in [C]) the injectivity radius of $\mathcal{S}_\lambda$ must be bounded above. This implies that the maximum of the sectional curvature of $\mathcal{S}_\lambda$ is bounded away from zero, which contradicts the fact that the scalar curvature $R$ approaches zero at infinity and the second fundamental form $II_\lambda$ is controlled by the scalar curvature. The theorem follows. ∎

Let $\mathcal{S}_\lambda(p)$ be the level set $\mathcal{S}_\lambda$ with $p \in \mathcal{S}_\lambda$ and $\mathcal{D}_\lambda(p)$ be the diameter of $\mathcal{S}_\lambda(p)$. By Gromov's compactness theorem, given a divergent sequence of points $\{p_i\}$ there is a subsequence $\{p_j\}$ such that

$$(\mathcal{M}^3, p_j, g/\mathcal{D}_\lambda^2(p_j)) \to \mathcal{P}_\infty$$

for some complete length space $\mathcal{P}_\infty$. Note that the length space $\mathcal{P}_\infty$ need not to be unique in general. It is known that $\mathcal{M}$ contains at least one ray $\gamma$ emanating from $\mathcal{O}$. The conclusion (4.2) implies that $\gamma$ determines a line in $\mathcal{P}_\infty$. By Corollary 7.10 in [Y], $\mathcal{P}_\infty$ splits isometrically,

$$\mathcal{P}_\infty = \mathbb{R} \times_1 \mathbb{Y}$$

where $\mathbb{Y}$ is some compact length space with $\operatorname{diam}(\mathbb{Y}) = 1$.

**Theorem 4.3** *Suppose that $(\mathcal{M}^3, g)$ is a gradient Ricci soliton with positive sectional curvature and that the scalar curvature assumes its maximum at an origin $\mathcal{O}$. Assume that $\mathbb{R} \times \Sigma^2$ cannot occur as a limit of dimension reduction. Then we have*

$$\limsup_{d(\mathcal{O},p) \to \infty} R\mathcal{D}_\lambda^2 < \infty. \tag{4.3}$$

**Proof.** Suppose this were not true. By Hamilton's dimension reduction, a non-flat factor is either (i) Hamilton's cigar $\Sigma^2$ or (ii) the standard 2-sphere $\mathbb{S}^2$. We shall claim that the 2-sphere $\mathbb{S}^2$ cannot occur as a non-flat factor, hence we get a contradiction. To do so, we choose a sequence of reference points and scale the metric as follows.

**Lemma 4.4** *If* $\limsup_{\lambda \to \infty} R\mathcal{D}_\lambda^2 = \infty$, *then there exists a sequence of points* $q_j \to \infty$, *a sequence of radii* $r_j$, *and a sequence of numbers* $\delta_j \to 0$ *such that*

(0) $\lim_{i \to \infty} R(q_i)\mathcal{D}_\lambda^2(q_i) = \infty$;

(a) $R(p) \leq (1 + \delta_j)R(q_j)$ *for all $p$ in the ball $B(q_j, r_j)$ of radius $r_j$ around $q_j$;*

(b) $r_j^2 R(q_j) \to \infty$;

(c) *if $s_j = d(q_j, \mathcal{O})$ is the distance of $q_j$ from some origin $\mathcal{O}$, then $\lambda_j = s_j/r_j \to \infty$;*

(d) *the balls $B(q_j, r_j)$ are disjoint.*

**Proof.** Pick a sequence $\varepsilon_j \to 0$, then choose $A_j \to \infty$ so that $A_j \varepsilon_j^2 \to \infty$ as well. As Lemma 22.2 in [H], let $\sigma_j$ be the largest number such that

$$\sup\{R(q)\mathcal{D}_\lambda^2(q) \mid \mathcal{D}_\lambda(q) \leq \sigma_j\} \leq A_j.$$

Since $R(q)\mathcal{D}_\lambda^2(q) \leq A_j$ if $\mathcal{D}_\lambda(q) \leq \sigma_j$, there exists some $q_j$ with

$$R(q_j)\mathcal{D}_\lambda^2(q_j) = A_j, \quad \mathcal{D}_\lambda(q_j) = \sigma_j, \quad \text{and } R(q_j) = \max_{q \in \mathcal{S}_{f(q_j)}} R(q),$$

otherwise $\sigma_j$ would not be maximal. Note that $\sigma_j \to \infty$ since $A_j \to \infty$ and $R$ is bounded.

Define $r_j = \varepsilon_j \sigma_j$. Conclusions (b) and (c) follow from the fact that $A_j \varepsilon_j^2 \to \infty$, $\varepsilon_j \to 0$ and $s_j \geq \sigma_j$ for $j$ sufficiently large. Next we check conclusion (a). If $p$ is in the ball of radius $r_j$ around $q_j$, either $\mathcal{D}_\lambda(p) \geq \sigma_j$ or $\mathcal{D}_\lambda(p) \leq \sigma_j$. In the first case, recall that $R$ decreases along any integral curve $\Gamma$ of the vector field $\nabla f/|\nabla f|^2$. We have from the choice of $q_j$

$$R(p) \leq \max_{q \in \mathcal{S}_\lambda(q_j)} R(q) = R(q_j)$$

which satisfies conclusion (a) with $\delta_j = 0$. In the second case, we need a lower bound for $\mathcal{D}_\lambda(p)$. Suppose that $a > b$ and consider for any point $q \in \mathcal{S}_a$ a minimizing geodesic $\gamma_q$ from $q$ to $\mathcal{S}_{b-1}$. Define $\beta_b = \max_{x \in \mathcal{S}_b} d(x, \mathcal{S}_{b-1})$, which is finite by compactness. By convexity of level sets, there are unique points $x_1 = $ the point of intersection of $\gamma_q$ with $\mathcal{S}_b$ and $x_2 = $ the point of intersection of $\gamma_q$ with $\mathcal{S}_{b-1}$. By the minimizing property of $\gamma_q$, the portion of $\gamma_q$ between $x_1$ and $x_2$ has length at most $\beta_b$. Note that $\beta_b \to 1$ as $b \to \infty$.

Again by convexity, the distance along $\gamma_q$ from $q$ to $x_1$ is less than or equal to $(a-b)\beta_b$. Thus $q$ has distance at most $\beta_b(a-b)$ from a point of $\mathcal{S}_b$ and by the triangle inequality

$$\mathcal{D}_b \geq \mathcal{D}_a - 2\beta_b(a-b) > \mathcal{D}_a - 3d(\mathcal{S}_a, \mathcal{S}_b). \tag{4.4}$$

We have from the choice of $\sigma_j$

$$R(p) \leq A_j/\mathcal{D}_\lambda^2(p)$$

and by (4.4)

$$\mathcal{D}_\lambda(p) \geq \mathcal{D}_\lambda(q_j) - 3r_j \geq (1-3\varepsilon_j)\sigma_j$$

so

$$R(p) \leq \frac{1}{(1-3\varepsilon_j)^2} \cdot \frac{A_j}{\sigma_j^2} = \frac{1}{(1-3\varepsilon_j)^2} R(q_j),$$

which satisfies conclusion (a) with $\delta_j = (1-3\varepsilon_j)^{-2} - 1$ and in either case $\delta_j \to 0$ as $\varepsilon_j \to 0$. Note that (a), (b) and (c) continue to hold if we pass to a subsequence. Any point $p$ in $B(q_j, r_j)$, the distance $d(p, \mathcal{O}) \geq s_j - r_j \to \infty$ as $j \to \infty$. Thus any fixed compact set does not meet the balls $B(q_j, r_j)$ for $j$ large enough. If we pass to a subsequence, the balls will all avoid each other. This proves the lemma. ∎

We scale the pointed Riemannian manifolds $(\mathcal{M}^3, q_i, g)$ both by $\mathcal{D}_\lambda^{-2}(q_i)$ and by $R(q_i)$. By hypothesis the cigar soliton cannot arise as a limit, thus, any non-flat factor must be the standard 2-sphere $\mathbb{S}^2$. Therefore, there is a subsequence $\{q_j\}$ such that the following sequences of pointed Riemannian manifolds

$$(\mathcal{M}^3, q_j, R(q_j)g/\pi^2) \text{ and } (\mathcal{M}^3, q_j, g/\mathcal{D}_\lambda^2(q_j))$$

converge to $\mathbb{R} \times_{1/\pi^2} \mathbb{S}^2$ and $\mathbb{R} \times_1 \mathbb{Y}$ respectively. On the other hand these two limit manifolds must be homothetic since both $\mathbb{S}^2$ and $\mathbb{Y}$ have finite diameters. We see that $\mathbb{Y}$ must be a 2-sphere of radius $1/\pi$, i.e., $\mathbb{Y} = \mathbb{S}^2_{1/\pi}$. This contradicts the choice of the reference points $q_i$ satisfying

$$\lim_{i \to \infty} R(q_i)\mathcal{D}_\lambda^2(q_i) = \infty,$$

and thus the 2-sphere $\mathbb{S}^2$ cannot form as a non-flat factor from dimension reduction as claimed. Therefore, the theorem follows. ∎

**Corollary 4.5** *Under the same hypothesis any sequence of pointed manifolds $(\mathcal{M}^3, p_i, g/\mathcal{D}_\lambda^2(p_i))$ always converges to $\mathbb{R} \times \mathbb{S}^2$ as long as $p_i$ goes to infinity and we have*
$$\limsup_{d(\mathcal{O},p) \to \infty} R\mathcal{D}_\lambda^2 = \pi^2.$$
*Moreover, there exist constants $c > 1$ and $s_0 > 0$ such that*
$$\frac{1}{c}\sqrt{s} < \mathcal{D}_\lambda < c\sqrt{s} \ \ \text{and} \ \ \frac{1}{c} < Rs < c$$
*for $s \geq s_0$. Therefore, the soliton opens like a paraboloid and the scalar curvature falls off like $1/s$.*

**Proof.** For convenience, let $\lambda_i$ be the value of $f$ evaluated at $p_i$, $\mathcal{D}_i$ the diameter $\mathcal{D}_{\lambda_i}$, and $g_i$ the metric $g/\mathcal{D}_i^2$. Since $\lim_{\lambda \to \infty} \mathcal{D}_\lambda/\lambda = 0$, there exists a divergent subsequence, say $\{p_j\}$, such that the sequence of pointed Riemannian manifolds $(\mathcal{M}^3, p_j, g_j)$ converges to a product space $\mathbb{R} \times_1 \mathbb{Y}$ with some length space $\mathbb{Y}$ of diameter 1. We want to claim that $\mathbb{Y}$ is a smooth manifold of dimension 2. By Gromov's compactness theorem the sequence of pointed Riemannian manifolds $(\mathcal{S}_j, p_j, g_j)$ converges to a compact set $\mathcal{S}_\infty \subset [-2, 2] \times \mathbb{Y}$, where $\mathcal{S}_j$ represents the level set $\mathcal{S}_{\lambda_j}$. To see that the sequence $(\mathcal{S}_j, p_j, g_j)$ does not collapse, it suffices to show that the positively curved surfaces $\mathcal{S}_j$ have uniformly bounded curvature, thus, the injectivity radii are uniformly bounded away from zero. Since the scalar curvature $R$ goes to zero at infinity, we have that

$$\begin{aligned} 0 < \mathcal{K}_\lambda &= K_\mathcal{M}(e_1, e_2) + \det(II_\lambda) \\ &= K_\mathcal{M}(e_1, e_2) + (Ric(e_1)Ric(e_2) - (Ric(e_1, e_2))^2)/|\nabla f|^2 \\ &< K_\mathcal{M}(e_1, e_2) + Ric(e_1) + Ric(e_2) < 2R \end{aligned}$$

for $\lambda$ sufficiently large. From the facts that $R\mathcal{D}_\lambda^2$ is bounded and $\mathcal{K}_\lambda < 2R$ for $\lambda$ sufficiently large, the sectional curvature of $\mathcal{S}_j$ with respect to $g_j$ is uniformly bounded. Therefore, the sequence $(\mathcal{S}_j, p_j, g_j)$ does not collapse, and the limit spaces $\mathcal{S}_\infty$ and $\mathbb{Y}$ are of dimension 2. For any constants $K, \delta > 0$, there exists some constant $\bar{\lambda} = \bar{\lambda}(K, \delta) > 0$ such that

$$(1 - \delta)d(\mathcal{S}_j, \mathcal{S}_\lambda) < |\lambda - \lambda_j| < d(\mathcal{S}_j, \mathcal{S}_\lambda)$$

for $\lambda_j > \bar{\lambda}$ and $\lambda \in (\lambda_j - K\mathcal{D}_j, \lambda_j + K\mathcal{D}_j)$. It follows that we have a uniform upper bound on the diameter $\mathcal{D}_\lambda$, $\lambda \in (\lambda_j - K\mathcal{D}_j, \lambda_j + K\mathcal{D}_j)$, with respect to the metric $g_j$. By local derivative estimates, the diameter estimate, and the fact that $R\mathcal{D}_\lambda^2$ is bounded, the sequence $(\mathcal{M}^3, p_j, g_j)$ converges smoothly to $\mathbb{R} \times \mathbb{Y}$. It is easy to see that $\mathbb{Y}$ is a 2-sphere with inner diameter 1. The first part of the corollary follows from the fact that the quantity $R(p_j)\mathcal{D}_\lambda^2(p_j)$ is scale invariant.

As a consequence of the first part we see that the sequence $(\mathcal{S}_{\lambda_j}, g/\mathcal{D}_{\lambda_j}^2)$ converges to $\mathbb{S}_{1/\pi}^2$ as long as $\lambda_j \to \infty$. It is known that $\mathcal{A}_\lambda$ grows linearly in $s$, thus there exists some constant $c > 1$ such that

$$\frac{1}{c}\sqrt{s} < \mathcal{D}_\lambda < c\sqrt{s} \text{ and } \frac{1}{c} < Rs < c.$$

∎

**Acknowledgements**. The author would like to thank Bing Cheng, Ben Chow, David Glickenstein, Peng Lu and Lei Ni for their helpful discussions and their interest in this work. The author also would like to thank the Institute for Advanced Study, UCLA, and the Center for Theoretical Sciences in Taiwan for the wonderful research environments.

# References

[BO]   R. Bishop and B. O'Neill, Manifolds with negative curvature, Trans. Amer. Math. Soc. Vol. **145**, 1-49.

[C]   C. Croke, Some isoperimetric inequalities and eigenvalue estimates, Ann. Sci. Éc. Norm. Sup. $4^c$ série, t. **13** (1980), 419-435.

[GLP]   M. Gromov, Structures métriques pour les variétés riemanniennes, edited by J. Lafontaine and P. Pansu, CEDIC, Paris, 1981

[GP]   K. Grove and P. Pertersen, Manifolds near the boundary of existence, J. Diff. Geom. **33** (1991), 379-394.

[GS]   R. Greene and K. Shiohama, Convex functions on complete noncompact manifolds: Topological structure, Invent. math. **63** (1981), 129-157.

[H]   R. S. Hamilton, The formation of singularities in the Ricci flow, Surveys in Differential Geometry **2** (1995), International Press, 7-136.

[I]   T. Ivey, New examples of complete Ricci solitons, Proc. Amer. Math. Soc. **122** (1994), no. 1, 241-245.

[K]   A. Kasue, A compactification of a manifold with asymptotically non-negative curvature, Ann. Sci. Ec. Norm. Sup., Paris, **21** (1988), 593-622.

[P]   G. Perelman, The entropy formula for the Ricci flow and its geometric applications, arXiv:math.DG/0211159.

[Y]   T. Yamaguchi, Collapsing and pinching under a lower curvature bound, Ann. of Math. **133** (1991), 317-357.